\documentclass[12pt]{amsart}
\usepackage{amsmath}

\usepackage{latexsym}
\usepackage{amssymb}

\usepackage[PostScript=dvips,noPS]{diagrams}

\def\square{\vbox{
      \hrule height 0.4pt
      \hbox{\vrule width 0.4pt height 5.5pt \kern 5.5pt \vrule width 0.4pt}
      \hrule height 0.4pt}}
\def\qed{{\unskip\nobreak\hfil\penalty50
  \hskip2em\hbox{}\nobreak\hfil\square
  \parfillskip=0pt\finalhyphendemerits=0\par\medbreak}}

\def\susp{\Sigma}

\renewcommand{\setminus}{\smallsetminus}

\begin{document}

\title{The space of commuting $n$-tuples in $SU(2)$}
\author{Thomas Baird}
\address{Mathematics Institute, University of Oxford\\UK}
\email{baird@maths.ox.ac.uk}
\urladdr{}
\author{Lisa C. Jeffrey}
\address{Department of Mathematics
University of Toronto \\ Canada}
\email{jeffrey@math.toronto.edu}
\author{Paul Selick}
\address{Department of Mathematics
University of Toronto \\ Canada}
\email{selick@math.toronto.edu}

 \newtheorem{theorem}{Theorem}[section]
 \newtheorem{corollary}[theorem]{Corollary}
 \newtheorem{conj}[theorem]{Conjecture}
 \newtheorem{prop}[theorem]{Proposition}
 \newtheorem{lemma}[theorem]{Lemma}
 \newtheorem{remark}[theorem]{Remark}
 \newtheorem{axiom}[theorem]{Axiom}
 \newtheorem{defn}[theorem]{Definition}
 \newtheorem{examp}[theorem]{Example}
 \newtheorem{proposition}[theorem]{Proposition}
 \newtheorem{claim}[theorem]{Claim}

\newcommand{\Proof}{{\noindent{\bf Proof:}\hspace{5pt}}}
\newcommand{\namedcong}[1]{\buildrel {#1} \over \cong}
\newcommand{\namedeq}[1]{\buildrel {#1} \over =}
\newcommand{\namedsimeq}[1]{\buildrel {#1} \over \simeq}
\newcommand{\homotopyequiv}{\simeq}
\newcommand{\heq}{\simeq}
\newcommand{\mapdown}[1]{\Big\downarrow\rlap{$\vcenter{\hbox{$\scriptstyle{#1}$}}$}}

\renewcommand{\implies}{\Rightarrow}
\newcommand{\assocto}{\longmapsto}
\newcommand{\beq}{\begin{equation}}
\newcommand{\eeq}{\end{equation}}
\renewcommand{\iff}{\Leftrightarrow}
\newcommand{\converges}{\rightarrow}
\newcommand{\contradiction}{\Rightarrow\Leftarrow}
\newcommand{\RR}{{\mathbb R}}
\newcommand{\CC}{{\mathbb C}}
\newcommand{\ZZ}{{\mathbb Z}}
\newcommand{\NN}{{\mathbb N}}
\newcommand{\QQ}{{\mathbb Q}}

 \newcommand{\reals}{\RR}
\newcommand{\complexes}{{\mathbb C}}
\newcommand{\integers}{{\mathbb Z}}
\renewcommand{\Proof}{{\noindent{\bf Proof:} }\hspace{5pt}}

\newcommand{\bprop}{\begin{proposition}}
\newcommand{\eprop}{\end{proposition}}
\newcommand{\bdef}{\begin{defn}}
\newcommand{\eedef}{\end{defn}}

\newcommand{\norm}[1]{||#1||}
\newcommand{\rp}{{\RR P}}
\newcommand{\circleaccent}[1]{\stackrel{\circ}{#1}}
\newcommand{\Int}{\mathop{\rm Int}\nolimits}
\newcommand{\eqvrel}{\mathord{\sim}}
\newcommand{\pt}{\mathord{*}}
\newcommand{\unique}{!}
\newcommand{\diam}{\mathop{\rm diam}\nolimits}
\newcommand{\im}{\mathop{\rm Im}\nolimits}
\newcommand{\Init}{\mathop{\rm Init}\nolimits}
\newcommand{\card}{\mathop{\rm card}\nolimits}
\newcommand{\coker}{\mathop{\rm coker}\nolimits}
\newcommand{\Hom}{\mathop{\rm Hom}\nolimits}
\newcommand{\sd}{\mathop{\rm sd}\nolimits}
\newcommand{\gsd}{\mathop{\rm SD}\nolimits}
\newcommand{\cell}{\mathop{\rm cell}\nolimits}
\newcommand{\id}{\mathop{\rm id}\nolimits}
\newcommand{\kk}{{\bf k}}
\newcommand{\degr}{\mathop{\rm deg}\nolimits}
\newcommand{\Obj}{\mathop{\rm Obj}\nolimits}
\newcommand{\rel}{\mathop{\rm rel}\nolimits}
\newcommand{\Map}{\mathop{\rm Map}\nolimits}
\newcommand{\hteq}{\simeq}
\newcommand{\smsh}{\wedge}

\newarrow{Equalto}{=}{=}{=}{=}{=}
\newarrowtail{Mono}{>}{<}{\vee}{\wedge}
\newarrow{Monoto}{Mono}{-}{-}{-}{>}

\begin{abstract}
Let $Y := \Hom(\ZZ^n, SU(2))$ denote the space of commuting $n$-tuples in $SU(2)$. We determine the homotopy type of the suspension $\susp Y$, and compute the integral cohomology groups of $Y$ for all positive integers $n$.
\end{abstract}
\maketitle

\section{Introduction}

It is interesting to study representations of discrete groups into compact Lie groups. For example, if $X$ is a smooth manifold then the space of homomorphisms $\Hom(\pi_1(X), G)$ may be identified with the space of `>flat $G$-connections on $X$ modulo based gauge transformations, which has diverse applications in geometry. 

An interesting special case is $X =(S^1)^n$, which has fundamental group $\ZZ^n$. 
The space of homomorphisms $\Hom(\ZZ^n, G)$ is identified with 
$$Y_G[n] := \{ (g_1, \dots, g_n) \in G^n \mid g_i g_j = g_j g_i 
\mbox{ for all $i,j$} \}.$$
(Usually $n$ and $G$ will be understood and omitted.)

This space was first studied by Adem and Cohen \cite{AC} (and has been further investigated with collaborators Torres-Giese \cite{ACT} and Gomez \cite{ACG} in connection with canonical filtrations of the classifying space $BG$).
They considered the problem in greater generality and obtained results even in
the more complicated case where $G$ is not compact.
 One of their main results (\cite{AC} Theorem 1.6) is a decomposition formula for the suspension of $Y_G[n]$ 
\begin{equation}\label{adcoh}
\susp Y_G[n] \simeq  \susp\left(\bigvee_{k=1}^n{n\choose k} Y_G[k]/S_k(G)\right) 
\end{equation}
where $S_k(G) \subset Y_G[k]$ consists of those k-tuples with some entry equal to the identity.
Although their papers focussed on the many interesting aspects of these spaces
rather than explicit cohomology computations, they also used their methods to
explicitly work out the cohomology groups of $Y_G[n]$ for $G=SU(2)$ in the
cases when $n=2$ and $n=3$.

Later, the first author \cite{Ba1}
gave a concise description of the 
rational cohomology ring as a ring of invariants under the action of the
Weyl group, and explicitly worked out Poincar\'e polynomials for the cases
$G=SU(2)$ and $SU(3)$. In the case $G=SU(2)$ and $n=2$ or $3$, the 
cohomology computed by~\cite{AC} agrees rationally
with that given by the first author.
The purpose of this paper is to calculate the stable homotopy
type of  $Y_{SU(2)}[n]$ for all $n$. Our formula is expressed in the form of (\ref{adcoh}), in terms explicit enough to compute (co)homology groups.
Our answer agrees rationally with that of the first author \cite{Ba1}.
It also agrees with \cite{AC} in the case $n=2$.
However for $n=3$ our results disagree with the published version of~\cite{AC} although they agree
rationally.
The authors would like to thank Alejandro~Adem and Fred~Cohen for the large
volume of email discussion during the interval since March~2006 when we first
emailed them our results pointing out the conflict with their paper.
Now that the issue has been resolved, we are pleased to publish our paper.

Since posting, we have been informed of an erratum~\cite{ACerr} 
to~\cite{AC}, and
of a 2008~preprint by M.C.~Crabb~\cite{Crabb} giving  results similar to ours.

We do not obtain the ring structure, although some information concerning the 
multiplication can be deduced from the rational calculation in \cite{Ba1}.  

\section{Wedge decomposition}
Let $G=SU(2)$ and let $G$ act on itself by conjugation.

Let
$T=\left\{\begin{pmatrix}\lambda&0\cr0&\lambda^{-1}\cr\end{pmatrix}\mid
\lambda\in S^1\right\}\subset G$ be a maximal torus.
The coset space $G/T=\{gT\}$ is homeomorphic to~$\CC P^1\cong S^2$.

Let $w=\begin{pmatrix}0&1\cr-1&0\cr\end{pmatrix}$.
Then $W=\{T=\bar{e},wT=\bar{w}\}$ is the Weyl group of~$G$.
The action of $W$ on $T$ is given by $\bar{w}\cdot t :=wtw^{-1}=t^{-1}$.
There is also an action of $W$ on $G/T$ given by
$\bar{w}\cdot \bar{g} :=\overline{wgw^{-1}}$, corresponding to the antipodal
action of $\ZZ/2$ on~$S^2$.

We use the subscripts $(~)_r$ and $(~)_s$ for the regular and singular
subsets.

Using this convention 
we set $T_s:=Z(G)=\{e,-e\}\subset T=T^W$ and $T_r:=T\setminus T_s$.
The action of $W$ on $T$ restricts to a free action on $T_r$ and the
trivial action on the fixed point set~$T_s=T^W$.
Similarly set
$$Y_s:=\{y\in Y\mid y_j\in Z(G)\ \forall j\}=(T^n)_s$$
and $Y_r:=Y\setminus Y_s$.
Note that $(T_s)^n = (T^n)_s = Y_s = \{\pm e\}^n$ is a collection
of $2^n$ isolated points in~$Y$.

Any set of commuting elements in $G$ must lie in a common maximal torus \cite{Ba1} and all maximal tori are conjugate. Consequently, the map
$$\phi = \phi[n]:(G/T)\times T^n\to Y[n]$$
satisfying
$$\phi(\bar{g},t_1,\ldots,t_n):= 
g\cdot(t_1,\ldots,t_n)=(gt_1g^{-1},\ldots,gt_ng^{-1})$$
is a $G$-equivariant surjection. The principal orbit type of $Y_r$ is $G/T$, so it follows that the restriction~$\pi$ of $\phi$,
$$\pi: (G/T) \times (T^n)_r \rightarrow Y_r$$ is a covering map. Because  
conjugacy classes in $G$ intersect $T$ in a $W$ orbit, we deduce that $\pi$ is a Galois cover, with deck transformation group $W$ acting diagonally on the product $G/T \times (T^n)_r$.

We thus obtain a homeomorphism,
\begin{equation}\label{homeomorphismtype:Y_r}
Y_r \cong \bigl((G/T) \times (T^n)_r\bigr)/W.
\end{equation}

Another way to look at this homeomorphism is as follows.
As noted above, for any $y\in Y$ there exists $\bar{g}\in G/T$ such that
$y\in\phi(\bar{g},T^n)$.
If $y\in Y_r$ then the class of $\bar{g}$ in $(G/T)/W\cong\rp^2$ is uniquely
determined by~$y$ so there is a well defined map $q:Y_r\to (G/T)/W$.
The map $q$ is a fibration with fibre $F=q^{-1}(\pt)=(T_r)^n$.
The Weyl group $W$ acts diagonally on $T^n$, and the
inclusion $F=(T_r)^n\subset T^n$ induces the $W$-action on~$F$.

Taking the pullback of $q$ with the universal covering projection
$\pi:S^2\to\rp^2$ gives a fibration
$\tilde{q}:\tilde{Y}_r\to S^2$.
\begin{diagram}
&&W&\rEqualto&W\cr
&&\dTo&&\dTo\cr
F&\rTo& \tilde{Y}_r&\rTo^{\tilde{q}}&S^2\cr
\dEqualto&&\dTo^{\pi}&&\dTo^{\pi}\cr
F&\rTo& Y_r&\rTo^{q}&\rp^2\cr
\end{diagram}

The action of $W$ on $\tilde{Y}_r\subset Y_r\times (G/T)$ is given by
$\bar{w}\cdot (y,\bar{g})=(y,\bar{w}\cdot\bar{g})=(y,\overline{wgw^{-1}})$.
This pullback fibration is trivial with retraction $r:\tilde{Y}_r\to F$
given by $r(y)=(t_1,\ldots t_n)$ where for all~$j$, $y_j=gt_jg^{-1}$ with
$\bar{g}=\tilde{q}(y)\in G/T=S^2$.
Thus $(r,q):\tilde{Y}_r\to  F\times S^2$ is a homeomorphism.
If $y_j=gt_jg^{-1}$ then we also have $y_j=gw(w\cdot t_j)(gwg^{-1})$.
Hence if $r(y,g)=(t_1,\ldots,t_r)$, then 
$r\bigl(w\cdot(y,g)\bigr)=(w\cdot t_1,\ldots,w\cdot t_r)$.
Therefore $\tilde{Y}_r\cong F\times S^2$ is a $W$-equivariant homeomorphism,
and so we obtain~(\ref{homeomorphismtype:Y_r}).

Recall that the set $Y_s = \{\pm 1\}^n$  is a collection of  $2^n$ isolated
points.
The following proposition shows
 that each point has a contractible neighbourhood in $Y$.
\begin{proposition}
The inclusion $Y_s\to Y$ is an absolute neighbourhood deformation
retract pair.
\end{proposition}
\Proof 
For $t\in T$ and $\epsilon>0$, let
$B_\epsilon(e)=\exp(-\epsilon,\epsilon)\subset T$
be a small interval about $e$ in~$T$.
For arbitrary $t\in T$, set $B_\epsilon(t)=tB_\epsilon(e)$.
For $y\in Y$,
if $(x,t)\in\phi^{-1}(y)$ is any particular preimage of~$y$, 
the set
$$V_\epsilon(y)=\bigcup_{g\in G}
\bigl(gB_\epsilon(t_1)g^{-1},\ldots, gB_\epsilon(t_n),g^{-1}\bigr),$$
is independent of the choice of $(x,t)$
and forms an open neighbourhood of $y$ in~$Y$.
$$\bigcap_{\epsilon>0}V_\epsilon(y)=\{y\}\cup\{\bar{w}\cdot y\}=W\cdot y.$$
In particular, if $y\in Y_s$ then 
$$\bigcap_{\epsilon>0}V_\epsilon(y)=\{y\}.$$
Set $V_s=\bigcup_{y\in Y_s}Y_y$, an open subset of~$Y$.

For small~$\epsilon$, define $H:B_\epsilon(e)\times I\to B_\epsilon(e)$
by~$H(\exp(x),s)=\exp(sx)$.
For arbitrary $t\in T$, this induces  a contraction
$H_t:B_\epsilon(t)\times I\to B_\epsilon(t)$ given by $H_t(x,s)=tH(t^{-1}x,s)$
of $B_\epsilon(t)$ to~$\{t\}$.
Set $H':=H_{-e}:B_\epsilon(-e)\times I\to B_\epsilon(-e)$.
Notice that $H$ and $H'$ are $W$-equivariant.
Suppose $y$ belongs to~$Y_s$.
Then for all~$j$, $y_j=\pm e\in T$.
For $v\in V_\epsilon(y)$, write $v=g\cdot(x_1,\ldots,x_n)$,
where $x_j\in B_\epsilon(y_j)$.
Define $H_y(v,s)=g\cdot\bigl(H_1(x_1,s),\ldots, H_n(x_n,s)\bigr)$
where $H_j=H$ or $H'$ according to whether $y_j=e$ or~$-e$.
Since $H$ and $H'$ are $W$-equivariant, the result is independent of the
choice of $g$ and $x_1,\ldots, x_n$, so produces a well defined contraction
$H_y:V_\epsilon(y)\times I\to V_\epsilon(y)$ of $V_\epsilon(y)$ to~$\{y\}$.
Therefore $V_s\simeq Y_s$.
Thus we have shown that $V_s$ is a neighbourhood deformation retract
of $Y_s$.
\qed

For a locally compact Hausdorff space~$X$, let $X^+$ denote its one-point
compactification.
Notice that $Y$ is the pushout $Y=Y_r\cup_{Y_r\cap V_s} V_s$.
Therefore 
\begin{equation} \label{e1} Y/V_s\simeq Y_r/(Y_r\cap V_s)\simeq Y_r^+.\end{equation}

 From equation~\ref{homeomorphismtype:Y_r},
\begin{equation} \label{e2} Y_r^+=(F\times S^2)^+/W=
\left(\frac{F^+\times S^2}{*\times S^2}\right)
\Big/W. \end{equation}

Since
$$F=(T_r)^n=
T^n\setminus Y_s,$$
it follows that
\begin{equation} \label{e3} F^+=T^n/Y_s. \end{equation}
In general,
\begin{equation} \label{e4}
\frac{B\times Z}{A\times Z}\cong\frac{B/A\times Z}{\pt\times Z}. \end{equation}
and
\begin{equation} \label{eq8}
(A/B)/W = (A/W)/ (B/W). \end{equation}
Therefore from (\ref{e2}) and (\ref{e3}) using (\ref{e4}) and (\ref{eq8})
we have 
\begin{equation} \label{e5}
 Y/V_s \cong 
\left(\frac{T^n\times S^2}
{Y_s\times S^2}\right)/W =\frac{\bigl((T^n\times S^2)/(\pt\times S^2)
\bigr)/W}
{(Y_s\times S^2/\pt\times S^2)/W}
.
\end{equation}
and after suspending,

\begin{equation} \label{e5.5}
 \Sigma(Y/V_s) \cong \frac{\Sigma(T^n\times S^2/\pt\times S^2)/W}
{\Sigma(Y_s\times S^2/\pt\times S^2)/W}
.
\end{equation}

\begin{remark} There is an ambiguity in the notation 
$X/K$ -- this might mean either the quotient by the action of a group $K$,
or the topological quotient where the subspace $K$ is collapsed to a point.
Unfortunately both notations are standard.
In the above equations, the quotients by~$W$ are those of  group actions and 
the others are
quotients of spaces.
\end{remark}
Since $W$ acts trivially on~$Y_s$, 
\begin{equation}  \label{e6}\left((Y_s\times S^2)/(\pt\times S^2)\right)/W
=(Y_s\times \rp^2)/(\pt\times \rp^2) \end{equation}
and after suspending we get
\begin{equation} \label{e7} \susp\Bigl ({(Y_s\times \rp^2)/(\pt\times \rp^2)} \Bigr )
\simeq\susp \Bigl ( (Y_s\smsh \rp^2)\vee Y_s \Bigr )
=\susp\left(\bigvee_{2^n-1}\rp^2 \vee Y_s\right). \end{equation}
This is  the denominator of (\ref{e5.5}).

We also have
 \begin{equation} \label{e8}\susp\left(\bigl(T^n\times S^2)/(\pt\times S^2)\bigr)/W
\right)
\simeq\susp\left(
\bigvee_{k=1}^n {n\choose k}\bigl((S^k\times S^2)/(\pt\times S^2)\bigr)/W
\right)
 \end{equation}
where $w$ acts  by the antipodal map on~$S^2$ and by the product
of $n$ reflections on $T^n=(S^1)^n$.
The left hand side of (\ref{e8}) is the numerator of (\ref{e5.5}).
Since the quotient map $(S^1)^k\to S^k$ is compatible with the $W$-action,
the action of $w$ on~$S^k$ has degree~$(-1)^k$ so is homotopic to the
negative of the antipodal map.

Next we need to identify the right hand side of (\ref{e8}).
Given a vector bundle~$\xi$, let $D(\xi)$ and $S(\xi)$ be its disk and
sphere bundles, and let $T(\xi):=D(\xi)/S(\xi)$ denote its Thom space.

\begin{lemma} \label{lthom}
Let $W$ act on $S^k$ and $S^2$ as above. Then
$$(S^k\times S^2)/(\pt\times S^2)\bigr)/W\simeq T(kL),$$ where $T(\xi)$ denotes
the Thom space of the bundle~$\xi$, and $L$ is the canonical line bundle
over~$\rp^2$.

\end{lemma}

\Proof
As $W$ spaces,
$\frac{S^k\times S^2}{\pt\times S^2}=
\frac{D^k\times S^2}{\partial D^k\times S^2}=T(k\tilde{L})$, where
$\tilde{L}$, a trivial line bundle, is the pullback of $L$ to~$S^2$
with $w$ acting by reflection.
Therefore 
$\left(\frac{S^k\times S^2}{\pt\times S^2}\right)\Big/W=T(k\tilde{L})/W=T(kL)$.
\qed

Putting this all together gives
\begin{equation} \label{e9}
\susp(Y/V_s)
\simeq\susp\left(\frac{\bigvee_{k=1}^n{n\choose k}T(kL)}
{\bigvee_{2^n-1}\rp^2 \vee Y_s}\right)
\end{equation}
Equation (\ref{e9}) is produced from (\ref{e5.5}) with replacements
from (\ref{e6}), (\ref{e7}), (\ref{e8}) and Lemma \ref{lthom}.

However
\begin{align} \label{e9a} 
\frac{\bigvee_{k=1}^n{n\choose k}T(kL)}
{\bigvee_{2^n-1}\rp^2 \vee Y_s}
&=\bigvee_{k=1}^n{n\choose k}\frac{T(kL)}{\rp^2}\Big/ Y_s\cr
&\simeq
\bigvee_{k=1}^n{n\choose k}\frac{T(kL)}{D(kL)}\Big/ Y_s\cr
&\simeq
\bigvee_{k=1}^n{n\choose k}\Sigma S(kL)\Big/ Y_s.\cr
\end{align}
Therefore \begin{equation} \susp(Y/V_s)
\simeq \susp\left(\bigvee_{k=1}^n{n\choose k}\Sigma S(kL)\Big/ Y_s\right)
\end{equation}
and since $V_s\simeq Y_s$ we get
\begin{equation} \label{e10} \susp Y\simeq
\susp\left(\bigvee_{k=1}^n{n\choose k}\Sigma S(kL)\right) \end{equation}

\section{Description of $\susp S(kL)$}

The goal of this section is to prove the following

\begin{prop}\label{prop3}
Decompose $\RR^{k+3} \cong \RR^3 \oplus \RR^k$, inducing disjoint embeddings of $\rp^2$ and $\rp^{k-1} $ into $\rp^{k+2}$. Then we have a homeomorphism
$$ \susp S(kL) \cong \rp^2 \backslash \rp^{k+2} \slash \rp^{k-1}$$
where the notation means that we contract $\rp^2$ and $\rp^{k-1}$ to distinct points. 
\end{prop}

We begin with a Lemma.

\begin{lemma}\label{lemma434}
If $L$ is the canonical line bundle over $\rp^m$
then $T(kL)\cong \rp^{k+m}/\rp^{k-1}$.
\end{lemma}

\Proof
Let $W=\ZZ/2$ acting via the antipodal action on $S^m$.
Let \mbox{$\pi:S^m\to\rp^m$}
 be the quotient map, and let $\tilde{L}=\pi^!L$ be the
pullback of $L$ to~$S^m$.
\begin{diagram}
W&\rEqualto&W\cr
\dTo&&\dTo\cr
S(k\tilde{L})&\rTo&D(k\tilde{L})\cr
\dTo^\pi&&\dTo^\pi\cr
S(k L)&\rTo&D(k\tilde{L})\cr
\end{diagram}

$\tilde{L}$ is a trivial bundle, so $D(k\tilde{L})\simeq I^k\times S^m$ and
$S(k\tilde{L})\simeq \partial(I^k)\times S^m$.
Therefore
$$T(kL)=\frac{D(kL)}{S(kL)}=\frac{D(k\tilde{L})}{S(k\tilde{L})}\Big/W
=\frac{I^k\times S^m}{\partial(I^k)\times S^m}\Big/W. $$

\begin{claim}
	
\begin{equation}\label{propquick} S^{m+k} 
\setminus S^{k-1} \cong { \Int} (D^k) \times S^m \end{equation}\end{claim}

\Proof 

$$\displaylines{S^{m+k} \setminus S^{k-1}  = \{ {\bf x} 
=(x_1 , \dots, x_k)\in{\RR}^k, 
\hfill \cr \hfill
{\bf y} = (y_1, \dots, y_{m+1}) \in {\RR}^{m+1} 
\,\Big  {|} \,
 |{\bf x}|^2 + |{\bf y}|^2 = 1, ~~ |{\bf x}|<1 \}.\cr} $$

Define a homeomorphism
$$ f: S^{m+k} \setminus S^{k-1} \to { \Int} (D^k) \times S^m $$
by 
$$ f({\bf x},{\bf y}) = ({\bf x}, {\bf y}/|{\bf y}|). $$
The inverse of this homeomorphism is
$$ g: { \Int} (D^k) \times S^m \to S^{m+k} \setminus S^{k-1} $$
defined by 
$$ 
g({\bf x},{\bf w}) = ({\bf x}, {\bf w}\sqrt{1 - |{\bf x}|^2}  ). $$ \
\qed

Taking one point compactifications of (\ref{propquick}) followed by $W$ orbits gives

$$ \rp^{m+k}/\rp^{k-1} \cong \frac{S^{m+k}}{S^{k-1}}/W \cong \frac{I^k \times S^m}{\partial I^k \times S^m}/W. $$

\noindent
\qed

\Proof(of Proposition \ref{prop3})
We have $\susp S(kL) \cong T(kL) /\rp^2$ where $\rp^2$ embeds into $T(kL)$ as the zero section. By Lemma \ref{lemma434} we know $ T(kL) \cong \rp^{k+2}/\rp^{k-1}$.  Under the homeomorphism (\ref{propquick}), the embedding of $\rp^{k-1}$ is by the $y$-coordinates and the embedding of $\rp^2$ is by the $x$-coordinates, so $\susp S(kL) \cong \rp^2 \backslash \rp^{k+2} \slash \rp^{k-1}$.

\qed

For cohomology calculations, the following corollary is convenient.

\begin{corollary}
$$\susp S(kL) \simeq
\begin{cases}S^3&{\rm if\ }k =1;\cr S^2\vee(\rp^{4}/\rp^2)&{\rm if\ }k=2\cr
\susp \rp^2\vee(\rp^{k+2}/\rp^{k-1})&{\rm if\ }k>2.
\end{cases}
$$
\end{corollary}
This combined with
\begin{equation} \label{e11}
\susp Y[n]\simeq
\susp\left(\bigvee_{k=1}^n{n\choose k}\Sigma S(kL)\right) \end{equation}
completely describes $\susp Y[n]$ up to homotopy equivalence.

We record for convenience that 
$$\tilde{H}^*( \rp^{k+2}/\rp^{k-1}) = \begin{cases} \ZZ&{\rm if\ }*= k+1 - (-1)^k;\cr \ZZ/2&{\rm if\ }*= k +\frac{3+(-1)^k}{2}  \cr 0&{\rm otherwise\ }\end{cases}$$
and
$$\tilde{H}^*( \rp^{4}/\rp^{2}) = \begin{cases}\ZZ/2&{\rm if\ }*= 4 \cr 0&{\rm otherwise\ }. \end{cases}$$

As noted earlier, although our results agree with the published 
version 
of~\cite{AC} when~$n=2$, they disagree when~$n=3$.
For reference, the following is an explicit listing of the groups in this case.

$$H^j( Y[3]) = \begin{cases} \ZZ &\mbox{if $j = 0 $} \cr
                            0 &\mbox{if $j = 1 $} \cr
                   \ZZ \oplus \ZZ \oplus \ZZ &\mbox{if $j = 2 $} \cr
           \ZZ \oplus \ZZ \oplus \ZZ \oplus \ZZ/2 &\mbox{if $j = 3 $} \cr
           \ZZ/2 \oplus \ZZ/2 \oplus \ZZ/2 \oplus \ZZ/2 &\mbox{if $j = 4 $} \cr
         \ZZ &\mbox{if $j = 5 $} \cr
         0 &\mbox{if $j >5 $} \cr
\end{cases} 
$$

\end{document}